\documentclass[11pt]{article}
\usepackage{graphicx}
\usepackage{amsmath,latexsym,amsbsy,amssymb, color}
\usepackage{graphicx}
\def\bpm{\begin{pmatrix}}
\def\bbm{\begin{bmatrix}}
\def\ebm{\end{bmatrix}}
\def\epm{\end{pmatrix}}

\def\argmin{\hbox{arg}\!\min}
\def\ds{\displaystyle}
\def\forall{\hbox{for all}~}
\def\L{{\bf L}}

\def\bfv{{\bf v}}

\def\bfe{{\bf e}}

\def\bfv{{\bf v}}
\def\ve{\varepsilon}
\def\n{\noindent}

\def\B{{\cal B}}

\def\Z{{\cal Z}}

\def\R{{\mathbb R}} 
\def\rank{{\rm Rank}}

\def\vp{\varphi}

\def\vs{\vskip 2em}

\def\v{\vskip 1em}

\def\V{{\cal V}}
\def\O{{\cal O}}

\def\M{{\cal M}}

\def\G{{\cal G}}
\def\C{{\cal C}}

\def\ov{\overline}

\def\Tilde{\widetilde}

\def\bega{\begin{array}}
\def\enda{\end{array}}
\def\begi{\begin{itemize}}
\def\endi{\end{itemize}}

\def\bel{\begin{equation}\label}
\def\eeq{\end{equation}}
\def\sqr#1#2{\vbox{\hrule height .#2pt
\hbox{\vrule width .#2pt height #1pt \kern #1pt
\vrule width .#2pt}\hrule height .#2pt }}
\def\square{\sqr74}
\def\endproof{\hphantom{MM}\hfill\llap{$\square$}\goodbreak}

\newtheorem{theorem}{Theorem}[section]

\newtheorem{lemma}{Lemma}[section]

\newtheorem{remark}{Remark}[section]
\newtheorem{definition}{Definition}[section]
\newtheorem{example}{Example}[section]

\begin{document}
\title{\bf  Generic uniqueness and conjugate points  for optimal control problems}\vs
	\author{Alberto Bressan$^{(1)}$, Marco Mazzola$^{(2)}$, and Khai T. Nguyen$^{(3)}$\\
\\
		{\small $^{(1)}$ Department of Mathematics, Penn State University}\\
		{\small $^{(2)}$ Sorbonne Universit\'e and Universit\'e Paris Cit\'e, CNRS, IMJ-PRJ, F-75005 Paris, France,}\\
		{\small $^{(3)}$  Department of Mathematics, North Carolina State University.}
		\\
		\quad\\
		{\small e-mails:  axb62@psu.edu,~~ marco.mazzola@imj-prg.fr,~~ khai@math.ncsu.edu}}
\maketitle
\begin{abstract} 
The paper is concerned with an optimal control problem on $\R^n$, where the dynamics is 
linear w.r.t.~the control functions.  For a terminal cost $\psi$ in a ${\cal G}_\delta$ set of $\C^4(\R^n)$
(i.e., in a countable intersection of open dense subsets), two main results are proved.
Namely: the set $\Gamma_\psi\subset\R^n$ 
of conjugate points
is closed, with locally bounded $(n-2)$-dimensional  Hausdorff measure.
Moreover, the set of initial points $y\in \R^n\setminus\Gamma_\psi$, which admit
 two or more globally optimal trajectories, is contained in
the union of a locally finite family of embedded manifolds.   In particular, the value function
is continuously differentiable on an open, dense subset of $\R^n$.

\quad\\
		\quad\\
		{\footnotesize
		\medskip
		
		\n {\bf AMS Mathematics Subject Classification.} 49K05, 49L12.
		}
	\end{abstract}

\maketitle

\section{Introduction}
\label{sec:1}
\setcounter{equation}{0}

Consider an optimal control problem of the form
\bel{oc1}
\hbox{minimize:}\qquad J^{y}[u]~\doteq~\ds\int_0^T L\bigl( x(t), u(t)\bigr)\, dt + \psi\bigl(x(T)\bigr),
\eeq
where the minimum is sought among all measurable functions $u:[0,T]\to\R^m$ and  $t\mapsto x(t)\in \R^n$ is the solution to a Cauchy problem
with dynamics linear w.r.t.~the control:
\bel{oc2} \dot x(t)~=~f\bigl(x(t), u(t)\bigr)~\doteq~f_0(x(t))+\sum_{i=1}^mf_i(x(t))\, u_i(t),\eeq
and initial data 
\bel{oc3} x(0)\,=\, y.\eeq
Here and in the sequel, the upper dot denotes a derivative w.r.t.~time.
Aim of this paper is to investigate properties of the set of conjugate points and of optimal 
controls which are {\it generic}.   Namely, for any given $f$ and $ L$, we seek properties 
which are valid for all terminal costs $\psi$ in a ${\cal G}_\delta$ set, i.e., in the  intersection
of countably many open dense subsets of a  Banach space~$\C^k(\R^n)$.

We assume that $f,L$ satisfy the following hypotheses:

\begi
\item[{\bf (A1)}] {\it The vector fields $f_i$, $i=0,1,\ldots,m$, are  three times continuously differentiable and 
satisfy the sublinear growth condition
 \bel{sublin}
 \bigl|f_i(x)\bigr|~\leq~c_1\bigl(|x|+1\bigr)
 \eeq
 for some constant $c_1>0$ and all $ x\in\R^n$.
 \item [{\bf (A2)}] The running cost  $L:\R^n\times\R^m\mapsto\R$ is three times continuously differentiable and satisfies
 \bel{Lbig}\left\{
 \bega{rl}
 L(x,u)&\geq~c_2\bigl( |u|^2-1\bigr),\\[2mm]
 |L_x(x,u)|&\leq~\ell(|x|)\cdot (1+|u|^2),
 \enda\right.
\eeq
 for some constant $c_2>0$ and some continuous function $\ell$. 
 Moreover,  $L$ is uniformly convex w.r.t.~$u$. Namely, for some $\delta_L>0$, 
the $m\times m$ matrix of second derivatives w.r.t.~$u$ satisfies
\bel{u-convex-L}
L_{uu}(x,u)~>~\delta_L\cdot \mathbb{I}_m \qquad\qquad\forall x,u.
\eeq
Here $ \mathbb{I}_m$ denotes the $m\times m$ identity matrix
}
\endi


Under the assumptions {\bf (A1)-(A2)}, the associated Hamiltonian $H:\R^n\times\R^n\to\R$ defined by 
\bel{Ham0}
H(x,p)~\doteq~\min_{\omega\in \R^m} \Big\{ L(x,\omega) + p\cdot f(x,\omega)\Big\}
\eeq
is three times continuously differentiable and   concave w.r.t.~$p$. 

Moreover, the Pontryagin necessary conditions \cite{BPi, Cesari, Clarke, FR} take the form
\bel{PMP1}
\left\{
\bega{rl} \dot x&=~H_p(x,p)~=~f\bigl(x,  u(x,p)\bigr),\\[3mm]
\dot p&=~~-H_x(x,p)=~- p\cdot  f_x\bigl(x,  u(x,p)\bigr)- L_x\bigl(x,u(x,p)\bigr),\enda
\right.
\eeq
where the optimal control $u(x,p)$ is determined as the pointwise minimizer
\bel{PMP2} u(x,p)~=~\argmin_{\omega\in \R^m} \Big\{ L(x,\omega) + p\cdot f(x,\omega)\Big\}.\eeq
Thanks to the above assumptions,  the system of ODEs (\ref{PMP1}) has continuously differentiable right hand side. 
In particular,  for any $z\in\R^n$, given the  terminal conditions
\bel{PMP3} x(T)~=~z,\qquad p(T)~=~\nabla\psi(z),\eeq
there exists a unique backward  solution  $t\mapsto \bigl(x (t,z),\, p(t,z)\bigr)$, at least locally in time.
The corresponding control in (\ref{PMP2}) will be denoted by
\bel{uxp}t~\mapsto~ u(t,z)~\doteq~u\bigl( x(t,z), \, p(t,z)\bigr).\eeq
\begin{definition}\label{def:1}
Given an initial point $y\in \R^n$, we say that 
a control $u^*:[0,T]\mapsto \R^m$ is a {\bf weak local minimizer} of the cost functional
\bel{Jbx} J^{y}[u]~\doteq~\int_0^T L\bigl( x(t), u(t)\bigr)\, dt + \psi\bigl(x(T)\bigr),
\eeq
subject to
\bel{dyn2}\dot x~=~f(x, u),\qquad\qquad 
x(0)~=~y,\eeq
if there exists $\delta>0$ such that
$
J^{y}[u^*]\leq J^{y}[u]$ for every measurable control $u(\cdot)$ such that
$\|u-u^*\|_{\L^\infty} <\delta$.
\end{definition}

Consider again the maps 
\bel{xpu}z\,\mapsto\, x(\cdot, z),\qquad z\,\mapsto\, p(\cdot, z),\qquad z\,\mapsto\, u(\cdot, z),
\eeq as  in (\ref{uxp}), obtained by solving the backward Cauchy problem
(\ref{PMP1})--(\ref{PMP3}).
Motivated by the classical definition of conjugate points~\cite{Cesari, Clarke, FR}, 
following \cite{CaFr14, PAC} we shall adopt
\begin{definition}\label{def:2} For the optimization problem 
(\ref{Jbx})-(\ref{dyn2}) a point $y \in \R^n$ is a {\bf  conjugate
point} if there exists $\ov z\in\R^n$ such that
$y = x(0,\ov z)$, the control
$u(\cdot,\ov z)$ is a weak local minimizer of (\ref{Jbx})-(\ref{dyn2}), and  moreover
\bel{Conj-def}
\hbox{\rm rank}\bigl[x_z(0,\ov z)\bigr]~<~n\,.
\eeq
\end{definition}
Here and in the sequel, by $x_z(t,\bar z)$ we denote the $n\times n$ Jacobian matrix of partial derivatives of the 
map $z\mapsto x(t,z)$ at the point $\bar z$.
It is understood that the above definition requires that the maps (\ref{xpu}) are well defined for all $t\in [0,T]$, 
as $z$ varies in a neighborhood of $\bar z$.

A necessary condition for the optimality of a trajectory starting at a conjugate point was recently given in \cite{BMN}. Relying on this result, a generic property for the set of conjugate points was obtained for a classical problem in the Calculus of Variations. In the present paper, we consider the more general case of an optimal control problem of the form (\ref{oc1}) and move forward in the study of the structure of the set of conjugate points, from the point of view of generic theory. Namely, for a generic terminal cost $\psi$: (i) the set $\Gamma_\psi\subset\R^n$ of all conjugate points
$y = x(0, \bar z)$,  such that $t\mapsto x(t,\bar z)$  is a globally optimal trajectory,
  is closed, with locally bounded $(n-2)$-dimensional Hausdorff measure, and
(ii) restricted to the open set $\R^n\setminus  \Gamma_\psi$, the set of initial points  where
the optimization problem has multiple global solutions is locally contained in the union of finitely many 
embedded manifolds of dimension $n-1$.  
As a consequence, the value function $V(y)$ for the optimal control problem (\ref{oc1})-(\ref{oc2}) is
continuously differentiable on an open, dense subset of $\R^n$.

We recall that, since the value function is semiconcave, it satisfies well known regularity 
properties~\cite{CaFr13, PAC, CaSin, CS}.
In particular, by the general result in \cite{AAC} the set of singularities of the value function
is contained in a countable union of Lipschitz manifolds
of dimension $n-1$. However, in some cases this set can be everywhere dense.
On the other hand, the present analysis shows that, in the generic case, this set of singularities
is nowhere dense.

For additional properties of conjugate points and their relations with optimal 
trajectories, we refer to \cite{CaFr14,CaFr96, ZZ}.

The remainder of the paper is organized as follows.   In Section~\ref{sec:2} we give a simple example
showing that, in a backward solution of (\ref{PMP1})--(\ref{PMP3}), the $p$-component can blow
up at an intermediate time $\tau\in \,]0,T[\,$.  On the other hand, Lemma~\ref{l:21} shows that,
for solutions which are globally optimal, both components $(x,p)$ satisfy a uniform bound for all 
$t\in [0,T]$.  Relying on the characterization of conjugate points recently given in \cite{BMN},
Section~\ref{sec:3} contains an analysis of the set of conjugate points, for a generic terminal cost 
$\psi\in \C^4(\R^n)$. In Section~\ref{sec:4}
the main results are stated in Theorem~\ref{t:41}.

\section{Bounds on solutions of   the PMP}
\label{sec:2}
\setcounter{equation}{0}
Even in the setting considered at {\bf (A1)-(A2)}, solutions to the backward Cauchy problem (\ref{PMP1})--(\ref{PMP3}) 
may not be well defined for all $t\in [0,T]$. Indeed, the $p$-component can blow up in finite time, 
as the following simple example shows.
\begin{example}\label{e:21} {\rm
Consider the scalar optimization problem on the time interval $t\in [0,2]$.
\bel{ex1} \hbox{Minimize:}\qquad \int_0^2  {u^2(t)\over 2}\, dt + \psi\bigl(x(2)\bigr),\eeq
subject to
\bel{ex2} \dot x~=~1+xu\,,\qquad\qquad x(0)=y.\eeq
Here the Hamiltonian function is
\bel{ex3} H(x,p)~=~\min_{\omega\in\R} \left\{ p\bigl(1+x\omega\bigr) + {\omega^2\over 2}\right\}~=~p-{p^2 x^2\over 2}\,,
\eeq
and the optimal control is 
\bel{ex4} u(x,p)~=~-px.
\eeq
The Hamiltonian system (\ref{PMP1}), (\ref{PMP3}) takes the form
\bel{PMP5}
\left\{
\bega{rl} \dot x&\ds=~1-px^2,\\[2mm]
\dot p&\ds=~p^2 x,\enda
\right.\qquad\qquad \left\{
\bega{rl} x(2,z)&\ds=~z,\\[2mm]
p(2,z)&\ds=~\psi'(z).\enda
\right.
\eeq
To solve the original optimization problem, the terminal point $z$ should be chosen so that $x(0)=y$.

We observe that, since the Hamiltonian function $H(x,p)$ is constant along trajectories of (\ref{PMP5}), 
 setting $w\doteq p x$ one finds 
\bel{Hconst}
p(t,z) -{w^2(t)\over 2}~=~H(x(t,z),p(t,z)) ~=~H(x(2,z),p(2,z))~=~H\bigl( z, \, \psi'(z)\bigr).
\eeq
We now choose $\psi$ to be a smooth, positive function such that $\psi'(-1)=2$.  For example
\bel{psidef}\psi(z)~=~2e^{z+1}
.\eeq
When $z=-1$, this choice yields $H\bigl( z, \, \psi'(z)\bigr)=0$.
By (\ref{Hconst})  and (\ref{PMP5}) it thus follows
\bel{w} p(t)~=~{2\over x^2(t)}\,,\qquad\qquad \dot x(t) ~=~-1,
\eeq
\[x(t)~=~1-t,\qquad\qquad
p(t) ~=~{2\over (1-t)^2}\,.
\]

Notice that in this  backward solution the $p$-component blows up as $t\to 1+$.
Computing the corresponding  control one finds
$$u(t)~=~- p(t) x(t)~=~{-2\over 1-t}\,,\qquad\qquad \int_0^2 {u^2(t)\over 2}\, dt ~=~+\infty.$$
Therefore, it cannot be optimal.
}
\end{example}

\begin{figure}[ht]
\centerline{
\hbox{\includegraphics[width=9cm]{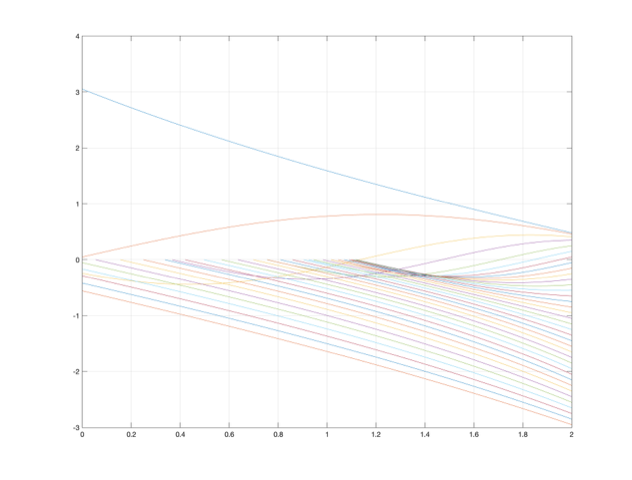}}  }
\caption{\small  Trajectories of (\ref{PMP5}), with terminal cost given by (\ref{psidef}).  
Notice that some backward trajectories terminate when  $x(t)=0$.  This happens when
the adjoint variable $p(t)$ blows up.   For $z>1/2$ backward solutions are still unique, but highly sensitive to the terminal data because the derivative $\psi'(z)= 2e^{z+1}$ has rapid growth. }
\label{f:ct3}
\end{figure}

While some of the solutions to the Hamiltonian system (\ref{PMP1})  can become unbounded in finite time, the following lemma shows
that the optimal ones satisfy uniform bounds.   A similar result was proved in \cite{BN}.

\begin{lemma}\label{l:21} Assume that  the couple $(f,L)$ satisfies {\bf (A1)-(A2)} and $\psi:\R^n\to \R_+$ is twice continuously differentiable. Then there exist continuous functions $\alpha,\beta,\gamma:\R_+\mapsto \R_+ $ such that  the following holds. Given any  initial point $y\in\R^n$ with $|y|\leq r$,
let  $u^*(\cdot)$ be a globally optimal control and let $x^*(\cdot)$  be the corresponding  optimal trajectory   for the 
problem (\ref{oc1})--(\ref{oc3}).   Then 
\bel{prio}\hbox{\rm ess-}\!\!\!\sup_{t\in [0,T]}  |u^*(t)|~\leq~\alpha(r),\qquad\qquad 
\sup_{t\in [0,T]}|x^*(t)|~\leq~\beta(r).
\eeq
Moreover, the adjoint variable $p^*(\cdot)$ satisfies
\bel{pbo}\sup_{t\in [0,T]}|p^*(t)|~\leq~\gamma(r).\eeq
\end{lemma}
{\bf Proof.} {\bf 1.} Consider any initial datum $y\in \ov B_r\doteq \bigl\{ x\in\R^n\,;~|x|\leq r\bigr\}$. 
Calling $x_0(\cdot)$  the solution  of (\ref{oc2}) with $u(t)\equiv 0$, i.e., the solution to
$$\dot x~=~f_0(x)\,,\qquad\qquad x(0)=y,$$
by (\ref{sublin}) it follows
\[
\sup_{t\in [0,T]}|x_0(t)|~\leq~(r+1)\cdot e^{c_1 T}.
\]
Let $(x^*,u^*)$ be a pair of optimal trajectory and optimal control of the optimization
problem (\ref{oc1})-(\ref{oc3}). By the first inequality in (\ref{Lbig}) it follows
\begin{multline}\label{u-L2}
\int_{0}^{T}\big|u^*(t)|^2dt~\leq~{1\over c_2}\,\left(\int_{0}^{T}L(x_0(t),0)dt+\psi(x_0(T))\right)+T\\
~\leq~{1\over c_2}\cdot\left(T\cdot \sup_{|x|\leq (r+1)\cdot e^{c_1T}} L(x,0)+\sup_{|x|\leq (r+1)\cdot e^{c_1T}}|\psi(x)|\right)+T~\doteq~\beta_1(r).
\end{multline}
Since   $x^*$ solves (\ref{oc2}) with $u\equiv u^*$, we have 
\[
|\dot{x}^*(t)|~\leq~c_1\cdot \bigl(|x^*(t)|+1\bigr)\cdot \left(1+\sum_{i=1}^{m}|u^*_i(t)|\right)~\leq~{c_1\over 2}\cdot 
\bigl(|x^*(t)|+1\bigr)\bigl(|u^*(t)|^2+m+2\bigr).
\]
Therefore, by (\ref{u-L2}) one obtains 
\[
\sup_{t\in [0,T]}|x^*(t)|~\leq~(r+1)\cdot \exp\left\{{c_1\over 2}\cdot \bigl[\beta_1(r)+(m+2)T\bigr]\right\}~\doteq~\beta(r).
\]

{\bf 2.} To derive a pointwise bound on $u^*$, for any constant $\lambda\geq 0$ we consider the truncated function
\[
u_{\lambda}(s)~=~\begin{cases}
u^*(s)&\mathrm{if}~~~~ s\notin   I_{\lambda},\\
0&\mathrm{if}~~~~ s\in I_{\lambda},
\end{cases}
\qquad\mathrm{with}\qquad I_{\lambda}\doteq \bigl\{s\in [0,T]: |u^*(s)|>\lambda\bigr\}.
\]
Calling $x_{\lambda}$  the solution of (\ref{oc2}) with $u\equiv u_{\lambda}$, we have 
\[
\sup_{t\in [0,T]}|x_{\lambda}(t)|~\leq~\beta(r),\qquad \sup_{t\in [0,T]}|x^*(t)-x_{\lambda}(t)|~\leq~\beta_2(r)\cdot \int_{I_{\lambda}}|u^*(s)|ds
\]
for some continuous function $\beta_2$. 
Recalling the constant $c_2$ in (\ref{Lbig}), for every 
\bel{labig}\lambda~\geq ~ c_2+\max_{|x|\leq \beta(r)} \bigl|L(x,0)\bigr|,\eeq 
we estimate the difference in the costs: 
$$\bega{l}\ds
0~\leq~J[u_{\lambda}]-J[u^*]~=~\int_{0}^{T}\Big[L\bigl(x_\lambda(t),u_\lambda(t)\bigr)-L\bigl(x^*(t),u^*(t)\bigr)\Big]\,dt+\psi\bigl(x_{\lambda}(T)\bigr)-\psi\bigl(x^*(T)\bigr)\\[4mm]
\qquad\ds
~\leq~\left(\bigl(T+\beta_1(r)\bigr)\cdot\sup_{|s|\leq \beta(r)}\ell(s)+\sup_{|x|\leq \beta(r)}\bigl|\nabla\psi(x)\bigr|\right)\cdot \beta_2(r)\cdot\int_{I_{\lambda}}|u^*(s)|ds\\[4mm]
\qquad\qquad \ds +\int_{I_{\lambda}}\Big[L(x^*(t),0)-L(x^*(t),u^*(t))\Big]dt\\[4mm]
\qquad\ds
~\leq~\alpha_1(r)\cdot  \int_{I_{\lambda}}\bigl|u^*(s)\bigr|\,ds-c_2\cdot \int_{I_{\lambda}}\bigl|u^*(s)\bigr|^2ds~\leq~\left(\alpha_1(r)-c_2\cdot\lambda\right)\cdot  \left(\int_{I_{\lambda}}\bigl|u^*(s)\bigr|\,ds\right)
\enda$$
for some continuous function $\alpha_1(\cdot)$.   

If now ess-sup$\,u^*(\cdot)>\lambda$, for some $\lambda$ which satisfies (\ref{labig}) together with
$\alpha_1(r)-c_2\cdot\lambda<0$, we obtain a contradiction.
We conclude that the first inequality in (\ref{prio}) must hold, with 
$$ \alpha(r)~=~\max\left\{{\alpha_1(r)\over c_2}\,,~ c_2+\max_{|x|\leq \beta(r)} \bigl|L(x,0)\bigr|\right\}.$$

{\bf 3.} Finally, a bound on the adjoint variable is obtained by observing that $p^*(\cdot)$ satisfies the linear ODE
\bel{linp} \dot p(t)~=~- p \cdot f_x\bigl( x^*(t), u^*(t)\bigr) - L_x\bigl( x^*(t), u^*(t)\bigr) ,\qquad\qquad p^*(T)~=~\nabla\psi( z).\eeq
By the previous steps, $x^*(t), u^*(t)$ remain uniformly bounded for $t\in [0,T]$, while the gradient 
$\nabla\psi(z)$ is bounded  as $z=x^*(T)$ ranges on the compact ball $\ov B_{\beta(r)}$. 
By standard estimates on solutions to linear systems of ODEs,  one obtains (\ref{pbo}).
\endproof

\begin{remark}\label{open}{\rm
It is easy to see that the set $\Omega$ of terminal points $z\in \R^n$ for which  the backward Cauchy problem (\ref{PMP1})--(\ref{PMP3}) is  well defined for all $t\in [0,T]$ is open.  Indeed, let $ z^*\in \Omega$, so that the corresponding solution $t\mapsto 
\bigl( x^*(t), p^*(t)\bigr)$ remains bounded:
\bel{sb1}\sup_{t\in [0,T]} \Big( \bigl|x^*(t)\bigr| + \bigl|p^*(t)\bigr|\Big) ~\doteq~R~<~+\infty.\eeq
Consider a smooth cutoff function $\vp:\R^n\times\R^n\mapsto [0,1]$ such that
$$\left\{\bega{rl} \vp(x,p)=1\quad &\hbox{if} ~~|x|+|p|\leq R+1,\\[1mm]
\vp(x,p)=0\quad &\hbox{if} ~~|x|+|p|\geq R+2. \enda\right.$$
Then the backward solutions to 
\bel{PPP}
\dot x\,=\,\vp(x,p)\,H_p(x,p),\qquad\qquad 
\dot p\,=\,-\vp(x,p)\,H_x(x,p),
\eeq
are well defined for all $t\in [0,T]$ and all terminal data (\ref{PMP3}).  In view of (\ref{sb1}),  by continuous dependence, 
for $|z-z^*|$ sufficiently small the solution to (\ref{PPP}) satisfies 
$$\sup_{t\in [0,T]} \Big( \bigl|x(t)\bigr| + \bigl|p(t)\bigr|\Big) ~<~R+1,$$
hence it is also a solution to the original system (\ref{PMP1}).
}
\end{remark}

\section{Generic properties of conjugate points in optimal control problems}
\label{sec:3}
\setcounter{equation}{0}
In this section we establish a generic property of the set of conjugate points for optimal control problems,
relying on the necessary condition stated in \cite[Theorem 2.2]{BMN}. 
More precisely, given a terminal cost 
$\psi\in\mathcal{C}^4(\R^n)$,  the conjugate points for the optimization problem (\ref{oc1})-(\ref{oc2}) are contained in the set 
$\big\{x(0,z)\,;~ (z,\bfv)\in\Omega_{\psi}\big\} $, with
\bel{J}  \Omega_{\psi}~\doteq~\Big\{(z,\bfv)\in\mathcal{O}_{\psi}\times S^{n-1}\,;~x_z(0, z)\bfv =0, \quad \bigl(p_z(0, z)\bfv\bigr)^{\dagger}\cdot x_{zz}(0,z)({\bf v}\otimes{\bf v})=0\Big\}.
\eeq
Here and in the sequel $S^{n-1}$ is the set of unit vectors in $\R^n$, 
the symbol $^\dagger$ denotes transposition, while $x_{zz}(0,z)$ denotes the second derivative
of the map $z\mapsto x(0,z)$ which is a symmetric, bilinear map from $\R^n\times \R^n$ into $\R^n$, and 
\[
\mathcal{O}_{\psi}~\doteq~\{z\in\R^n;~\mathrm{the~system}~(\ref{PMP1})-(\ref{PMP3})~\mathrm{admits~a~solution~on~} [0,T]\}.
\]
%
%
%
Notice that by Remark \ref{open}, $\mathcal{O}_{\psi}$ is an open set in $\R^n$. Our  main result of this section can be stated as follows.


\begin{theorem}\label{t:31}  Let $f_0,\ldots, f_m$ and $L$ be given, satisfying 
{\bf (A1)}-{\bf (A2)}. Then there exists a $\mathcal{G}_{\delta}$ subset $\mathcal{M}\subseteq\mathcal{C}^4(\R^n)$ such that, for every $\psi\in \mathcal{M}$, the set $\Omega_{\psi}$  in (\ref{J})
  is  
   an embedded manifold  of dimension $n-2$. 
  \end{theorem}
  
{\bf Proof.} {\bf 1.}  We shall first assume that  the Hamiltonian system  (\ref{PMP1})-(\ref{PMP3}) admits a global solution for all $z\in \R^n$ and $\psi\in \mathcal{C}^4(\R^n)$. For a given terminal cost  $\psi\in\mathcal{C}^4(\R^n)$, we define the $\C^1$ map  $\Phi^{\psi}:\R^d\times S^{n-1}\mapsto\R^n\times\R$ 
by setting
\bel{Phi-psi}
\Phi^{\psi}(z,{\bf v})~\doteq~\left(x_z(0,z){\bf v} \,, ~ \bigl(p_z(0, z)\bfv\bigr)^{\dagger}\cdot x_{zz}(0,z)({\bf v}\otimes{\bf v})\right).
\eeq
For $k\geq 1$, let $\overline{B}_k\subset\R^n$  denote the closed ball centered at the origin with radius $k$. We then define the subset of $\mathcal{C}^4(\R^n)$ as
\bel{M-k}
\mathcal{M}_{k}~\doteq~\left\{\psi\in \mathcal{C}^4(\R^n) \,;~~ \Phi^{\psi}\Big|_{\overline{B}_k\times S^{n-1}}\text{ is transversal to } \{\bf0\}\right\}.
\eeq
Here $\{\bf0\}$ denotes  the zero-dimensional manifold containing the single point $(0,0)\in \R^n\times\R$. Otherwise stated, $\M_k$ is the set of all terminal costs
$\psi\in \mathcal{C}^4(\R^n)$ with the following property.
Setting  $$\mathcal{Z}_0~\doteq~\bigl\{(z,{\bf v})\in \overline{B}_k\times S^{n-1}\,;~\Phi^\psi(z,\bfv)=\bf0\bigr\},$$
and denoting by $D\Phi^\psi$ the $(n+1)\times 2n$ Jacobian matrix of partial derivatives of the map
$(z, \bfv)\mapsto\Phi^\psi(z,\bfv)$,
 there holds
\bel{F-r1}
\mathrm{rank}\, \bigl[D\Phi^\psi(z,\bfv)\bigr]\,=\,n+1\qquad\qquad \forall (z,{\bf v})\in \mathcal{Z}_{0}\,.
\eeq
Note that the map $\Phi^\psi$ in (\ref{Phi-psi}) has polynomial dependence w.r.t.~$\bfv$. Hence it is 
well defined for all $\bfv\in\R^n$.
\v
{\bf 2.} We claim that each set $\M_k\subset\C^4(\R^n)$ is open. Indeed, 
 by the continuity of $D\Phi^\psi$, the identity (\ref{F-r1}) remains valid for every $(z,{\bf v})\in\O$, where $ \O$ is an open neighborhood of  $\mathcal{Z}_{0}$. 
By the definition of $\Z_0$ it follows
\bel{n-ze}
\min_{(z,{\bf v})\in \B_k\times S^{n-1}\backslash \O}
\big|\Phi^\psi(z,\bfv)\big|~>~0.
\eeq
Again, since the map $\psi\mapsto \big(\Phi^{\psi}, D\Phi^{\psi}\big)$ is continuous, there exists $\ve>0$ sufficiently small such that for every terminal cost 
$\tilde \psi$ such that $\big\|\tilde \psi-\psi\big\|_{\C^4}<\ve$, both (\ref{F-r1}) and (\ref{n-ze})
remain valid for $\tilde\psi$ as well. 
We thus conclude that the set $\M_k$ is open in $\C^4(\R^n)$.
%
%
%
%
\v
{\bf 3.} In the next steps we will prove that each set $\M_k$ is also dense in $\C^4(\R^n)$. 

For this purpose, fix any $\Tilde{\psi}\in \C^4(\R^n)$ and let $\ve>0$ be given.
We first approximate $\Tilde{\psi}$ by a smooth function $\psi\in\C^\infty$ such that 
$\big\|\psi-\Tilde\psi\big\|_{\C^4}<\ve$. 

Next, we construct a perturbed 
function $\psi^\theta$, arbitrarily close to $\psi$ in the $\C^4$ norm, which 
lies in $\M_k$.  Toward this goal, we first define $\vp:\R^{n^2+n}\times\R^n\to\R^n$ by setting
\bel{vp}
\vp(\theta,z)~\doteq~\sum_{i,j=1}^{n}\theta_{ij} \,z_iz_j+\sum_{k=1}^{n}\theta_k\,z_k^3,\qquad\qquad \theta=(\theta_{ij}, \theta_k)\in \R^{n^2+n}, ~~z\in \R^n.
\eeq
A straightforward differentiation  yields
\bel{d-vp}
{\partial\over \partial z_i}\vp(\theta,z)\,=\,\sum_{j=1}^n (\theta_{ij}+\theta_{ji})  z_j + 3  \theta_i z_i^2,\qquad {\partial^{2}\over \partial z_i\partial z_j}\vp(\theta,z)\,=\,\left\{\bega{cl} 2\theta_{ii}+ 6  \theta_i z_i\quad &\hbox{if} \quad i=j,\\[2mm]
\theta_{ij}+\theta_{ji}~~ &\hbox{if} \quad i\not= j.\enda\right.
\eeq
Fix a point $\bar z\in \B_k$, and
consider the family of perturbed functions defined by
\bel{psite}\psi^\theta(z)~\doteq~\psi(z) +\eta(z-\bar{z})\cdot\vp(\theta,z-\bar{z}),\qquad \qquad (\theta,z)\in \R^{n^2+n}\times\R^n,
\eeq
where $\eta:\R^n\mapsto [0,1]$ is a smooth cutoff function such that 
\bel{cutoff}\eta(y)~=~\left\{ \bega{rl} 1\quad &\hbox{if}\quad |y|\leq 1,\\[2mm]
0\quad &\hbox{if}\quad |y|\geq 2.\enda\right.\eeq
For any $(\theta,z)\in \R^{n^2+n}\times\R^n$, we denote by  $t\mapsto \bigl(x^\theta(t,z),p^\theta(t,z)\bigr)$ 
the solution of (\ref{PMP1}) with terminal conditions 
\bel{T-con}
x^{\theta}(T,z)\,=\,z,\qquad p^{\theta}(T,z)\,=\,\nabla\psi^\theta(z).
\eeq
\v
{\bf 4.}
We claim that the map 
$$ (z,\bfv, \theta)~\mapsto~\Phi^{\psi^\theta}(z,\bfv)~\doteq~ \Big(x^{\theta}_z(0,z){\bf v} , ~
\bigl(p^\theta_z(0, z)\bfv\bigr)^{\dagger}\cdot x^\theta_{zz}(0,z)({\bf v}\otimes{\bf v})\Big)$$
is transversal to $\{{\bf 0 }\}\subset\R^n\times\R$ at the point $(\bar z, \bfv, 0)$. This will certainly be true if  either $\Phi^{\psi^\theta}(\bar{z},\bfv)_{|_{\theta=0}}\not={\bf 0}$, or else
if the Jacobian matrix $D_\theta \Phi^{\psi^{\bar{z}}(\theta,\cdot)}(\bar z, \bfv)_{|_{\theta=0}}$
of partial derivatives w.r.t.~$\theta_{ij}, \theta_k$ has  rank $n+1$.
\medskip

Assume that $\Phi^{\psi^\theta}(z,\bfv)_{|_{\theta=0}}={\bf 0}$. Then 
    we have 
\bel{Miss}
x^{\theta}_z(0,\bar{z})_{|_{\theta=0}}{\bf v}~=~0.
\eeq

Differentiating the Pontryagin equations (\ref{PMP1}) and the terminal conditions 
(\ref{T-con})  w.r.t.~$z$ and then w.r.t.~$\theta_{ij}$, we obtain
\bel{Ham2}
\left\{\bega{rl} \dfrac {d}{dt}x_z^\theta(t, z)&=\, H_{xp}\cdot x_z^\theta(t, z)+H_{pp}\cdot p_z^\theta(t, z),\\[4mm]
\dfrac {d}{dt} p_z^\theta(t, z)&=-H_{xx}\cdot x_z^\theta(t, z)-H_{px}\cdot p_z^\theta(t, z),\enda\right.
\qquad
\left\{\bega{rl} x_z^\theta(T, z)&=~\mathbb{I}_n,\\[4mm]
p_z^\theta(T, z) &=\,\psi^\theta_{zz}(z),\enda\right.\eeq
\bel{Ham3}
\left\{\bega{rl} \dfrac {d}{dt} x_{\theta_{ij}}^\theta(t, z)&=\, H_{xp}\cdot x_{\theta_{ij}}^\theta(t, z)+H_{pp}\cdot p_{\theta_{ij}}^\theta(t, z),\\[4mm]
\dfrac {d}{dt} p_{\theta_{ij}}^\theta(t, z)&=- H_{xx}\cdot x_{\theta_{ij}}^\theta(t, z)-H_{px}\cdot p_{\theta_{ij}}^\theta(t, z),\enda\right.
\ 
\left\{\bega{rl} x_{\theta_{ij}}^\theta(T, z)&=~0,\\[4mm]
p_{\theta_{ij}}^\theta(T, z) &=~D_z\vp_{\theta_{ij}}(\theta,z-\bar{z}).\enda\right.\eeq
Since $p_{\theta_{ij}}^\theta(T, \bar{z})=\vp_{\theta_{ij}}(\theta,\bar{z})=0$, one has
\[
x_{\theta_{ij}}^\theta(t,\bar z)~=~p_{\theta_{ij}}^\theta(t,\bar z)~=~0,\qquad t\in [0,T].
\]

For any  ${\bf v}= (\bfv_1,\ldots, \bfv_n)\in S^{n-1}$, differentiating (\ref{Ham2}) w.r.t.~$\theta_{ij}$ we thus obtain
\bel{Ham4}
\left\{\bega{rl} \dfrac {d}{dt}\big( x_{z}^\theta(t,\bar z)\bfv\big)_{\theta_{ij}}&=\, H_{xp}\cdot \big( x_{z}^\theta(t,\bar z)\bfv\big)_{\theta_{ij}}+H_{pp}\cdot \big( p_{z}^\theta(t,\bar z)\bfv\big)_{\theta_{ij}},\\[4mm]
\dfrac {d}{dt}\big( p_{z}^\theta(t,\bar z)\bfv\big)_{\theta_{ij}}&= -H_{xx}\cdot \big( x_{z}^\theta(t,\bar z)\bfv\big)_{\theta_{ij}}-H_{px}\cdot \big( p_{z}^\theta(t,\bar z)\bfv\big)_{\theta_{ij}},\enda\right.
\eeq
\bel{T1}
\big( x_{z}^\theta(T,\bar z)\bfv\big)_{\theta_{ij}}~=~0,\qquad \big( p_{z}^\theta(T,\bar z)\bfv\big)_{\theta_{ij}} ~=~ {\bf v}_j\cdot {\bf e}_i+{\bf v}_i\cdot {\bf e}_j\,.
\eeq
As usual,  $\{{\bf e}_1,\cdots, {\bf e}_n\}$ denotes the standard basis of $\R^n$. 

For every $i\in \{1,\cdots,n\}$, let $S_i(t)$ and  $P_i(t)$ be the $n\times n$ matrices with columns $\left\{\big( x_{z}^\theta(t,\bar z)\bfv\big)_{\theta_{ij}}\right\}_{j=1,\ldots,n}$ and $\left\{\big( p_{z}^\theta(t,\bar z)\bfv\big)_{\theta_{ij}}\right\}_{j=1,\ldots,n}$, respectively.
That means:
\bel{SPi}
\begin{cases}
S_i(t)~=~\left(\partial_{\theta_{i1}} 
x_{z}^\theta(t,\bar z)_{|_{\theta=0}}\bfv\big)\Bigg|~\cdot~\cdot~\cdot~\Bigg|
\partial_{\theta_{in}} 
x_{z}^\theta(t,\bar z)_{|_{\theta=0}}\bfv\right),\\[6mm]
P_i(t)~=~\left(\partial_{\theta_{i1}} 
p_{z}^\theta(t,\bar z)_{|_{\theta=0}}\bfv\big)\Bigg|~\cdot~\cdot~\cdot~\Bigg|
\partial_{\theta_{in}} 
p_{z}^\theta(t,\bar z)_{|_{\theta=0}}\bfv\right),
\end{cases}
\eeq
for all $ t\in [0,T]$.
For notational convenience, we now introduce the functions
\[
A(t)\,\doteq\,H_{xx}\bigl(x(t,\bar{z}),p(t,\bar{z})\bigr),\qquad B(t)\,\doteq\,H_{xp}\bigl(x(t,\bar{z}),p(t,\bar{z})
\bigr),\qquad C(t)\,\doteq\,H_{pp}\bigl(x(t,\bar{z}),p(t,\bar{z})\bigr).
\]
By (\ref{Ham2}) and (\ref{Ham4})-(\ref{T1}), 
 for every $i\in \{1,\ldots,n\}$  the maps $t\mapsto \bigl(S_i(t),P_i(t)\bigr)$  and $t\mapsto \bigl(x_z^\theta(t, \bar{z})_{|_{\theta=0}}, ~p_z^\theta(t, \bar{z})_{|_{\theta=0}} \bigr)$   provide a solution to the system of ODEs
\bel{Ham-551}
\begin{cases}
 \dot{X}(t)&=\, B(t) X(t)+C(t) Y(t),\\[4mm]
\dot{Y}(t)&=\, -A(t)X(t)-B(t)Y(t),
\end{cases}
\eeq
with terminal data
\[
X(T)~=~0,\qquad\qquad Y(T)~=~\Big({\bf v}_1\, {\bf e}_i+{\bf v}_i\, {\bf e}_1,~~\cdots~~, {\bf v}_n\, {\bf e}_i+{\bf v}_i\, {\bf e}_n\Big),
\]
\[
X(T)~=~\mathbb{I}_n,\qquad\qquad Y(T)~=~D^2\psi(\bar{z}),
\]
 respectively.
Since  ${\bf v}\in S^{n-1}$ is a unit vector, for some index
$i_0\in\{1,\cdots, n\}$ one must have ${\bf v}_{i_0}\neq 0$. We claim that 
\bel{r-S}
\mathrm{rank}[S_{i_0}(0)]~=~n.
\eeq
Indeed, if the $n\times n$ matrix $S_{i_0}$ does not have full rank, there exists some ${\bf w}\in S^{n-1}$ 
such that $S_{i_0}(0){\bf w}=0$. To derive a contradiction, notice that the maps $t\mapsto \big(S_{i_0}(t){\bf w},P_{i_0}(t){\bf w}\big)$ and  $t\mapsto \big(x^{\theta}_z(t,\bar{z})_{|_{\theta=0}}{\bf v},p^{\theta}_z(t,\bar{z})_{|_{\theta=0}}{\bf v}\big)$ are the solutions of (\ref{Ham-551}) with terminal data
\bel{cd11}
X(T)~=~0,\qquad\qquad Y(T)~=~{\bf v}_{i_0} {\bf w}+{\bf v^{\dagger}{\bf w}}\cdot {\bf e}_{i_0}~\neq~0,
\eeq
\bel{cd12}
X(T)~=~{\bf v},\qquad\qquad Y(T)~=~D^2 \psi(\bar{z}){\bf v}.
\eeq
From (\ref{Miss}), (\ref{cd11}) and (\ref{cd12}), it now follows 
\[
P_{i_0}(0){\bf w}~\neq~0,\qquad \qquad p^{\theta}_z(0,\bar{z})_{|_{\theta=0}}{\bf v}~\neq~0.
\]
 In this case, there exists an invertible $n\times n$ matrix  $M$    such that 
\[
\big[ P_{i_0}(0){\bf w}\big] ~=~p^{\theta}_z(0,\bar{z})_{|_{\theta=0}}{\bf v}M~\neq~0.
\]
Therefore the two maps $$t\,\mapsto\, \bigl(S_{i_0}(t){\bf w},\,P_{i_0}(t){\bf w}\bigr),\qquad\qquad 
t\,\mapsto\, \Big(x^{\theta}_z(t,\bar{z})_{|_{\theta=0}}{\bf v}M,~p^{\theta}_z(t,\bar{z})_{|_{\theta=0}}{\bf v}M\Big),$$
both solve  (\ref{Ham-551}) with the same initial data 
\[
X(0)~=~0,\qquad \qquad Y(0)~=~\big[ P_{i_0}(0){\bf w}\big] ~=~p^{\theta}_z(0,\bar{z})_{|_{\theta=0}}{\bf v}M.
\]
By uniqueness, we obtain 
\[
\big(S_{i_0}(t){\bf w},P_{i_0}(t){\bf w}\big)~=~\bigr(x^{\theta}_z(0,\bar{z})_{|_{\theta=0}}{\bf v}M~,~p^{\theta}_z(0,\bar{z})_{|_{\theta=0}}{\bf v}M\bigl)\qquad\forall t\in [0,T].
\]
In particular, from (\ref{cd11}) and (\ref{cd12}) it follows
\[
{\bf v}M~=~x^{\theta}_z(T,\bar{z})_{|_{\theta=0}}{\bf v}M~=~S_{i_0}(T){\bf w}~=~0,
\]
reaching a contradiction. Hence, (\ref{r-S}) holds.
\v
{\bf 5.} Next, differentiating (\ref{PMP1}) w.r.t.~$\theta_k$, one checks that the map 
$t\mapsto \bigl(x_{\theta_k}^\theta(t, \bar{z}),p_{\theta_k}^\theta(t, \bar{z})\bigr)$ provides  a solution of (\ref{Ham-551}) with terminal data
\[
\partial_{\theta_k}x^\theta(T, \bar{z})\,=\,0,\qquad \qquad \partial_{\theta_k}p^\theta(T, z)\,=\,\partial_{\theta_k}\vp(\theta,0)\,=\,0.
\]
This implies
\bel{zero1}
x_{\theta_k}^\theta(t,\bar z)~=~p_{\theta_k}^\theta(t,\bar z)~=~0\quad \qquad\forall t\in [0,T].
\eeq
Hence, differentiating (\ref{Ham2}) w.r.t.~$\theta_k$,  we check that  the map $t\mapsto \big(\partial_{\theta_k}x_{z}^\theta(t, \bar{z}),~\partial_{\theta_k}p_{z}^\theta(t, \bar{z})\big)$ provides  a solution of (\ref{Ham-551}) with terminal data
\[
\partial_{\theta_k}x_z^\theta(T, \bar{z})~=~0,\qquad \qquad \partial_{\theta_k}p_z^\theta(T, \bar{z})~=~\partial_{\theta_k}\vp_z(\theta,0)~=~0.
\]
Therefore
\bel{zero2}
\partial_{\theta_k}x_{z}^\theta(t,\bar z)~=~\partial_{\theta_k}p_{z}^\theta(t,\bar z)~=~0,\qquad\forall t\in [0,T].
\eeq
For $k\in \{1,\cdots,n\}$, forming an augmented $n\times (n+1)$ matrix we set
\bel{S-k}
\begin{split}
\Tilde{S}_{i_0,k}(0)&~\doteq~\left(\partial_{\theta_{i_01}} 
x_{z}^\theta(t,\bar z)_{|_{\theta=0}}\bfv\big)\bigg|\cdot\cdot\cdot\bigg|
\partial_{\theta_{i_0 n}} 
x_{z}^\theta(t,\bar z)_{|_{\theta=0}}\bfv \Bigg|~\partial_{\theta_k} 
x_{z}^\theta(t,\bar z)_{|_{\theta=0}}\bfv\right)\\
&~=~ \begin{bmatrix}~~S_{i_0}(0)&\Big|~
0_{n\times 1}\end{bmatrix}.
\end{split}
\eeq

Using (\ref{zero1}) and (\ref{zero2}), differentiating twice (\ref{Ham2}) w.r.t $\theta_k, z$, we obtain  
\bel{Ham7}
\begin{cases}
\dfrac {d}{dt}\big[x_{zz}^\theta(t,\bar z)(\bfv\otimes\bfv)\big]_{\theta_{k}}~=\,B(t)\cdot \big[ x_{zz}^\theta(t,\bar z)(\bfv\otimes\bfv)\big]_{\theta_{k}}+C(t)\cdot \big[ p_{zz}^\theta(t,\bar z)(\bfv\otimes\bfv)\big]_{\theta_{k}},\\[4mm]
\dfrac {d}{dt}\big[p_{zz}^\theta(t,\bar z)(\bfv\otimes\bfv)\big]_{\theta_{k}}~=-A(t)\cdot \big[x_{zz}^\theta(t,\bar z)(\bfv\otimes\bfv)\big]_{\theta_{k}}-B(t)\cdot \big[ p_{zz}^\theta(t,\bar z)(\bfv\otimes\bfv)\big]_{\theta_{k}},
\end{cases}
\eeq
with  terminal conditions
\bel{T-2}
\big[x_{zz}^\theta(T,\bar z)(\bfv\otimes\bfv)\big]_{\theta_{k}}~=~0,\qquad \qquad\big[p_{zz}^\theta(T,\bar z)(\bfv\otimes\bfv)\big]_{\theta_{k}}~=~6\bfv_k^2\bfe_k\,.
\eeq
By \eqref{zero2} one has the identities
\[\begin{cases}\big[\bigl(p^\theta_z(t, \bar z)\bfv\bigr)^{\dagger}\cdot x^\theta_{zz}(t,\bar z)({\bf v}\otimes{\bf v})\big]_{\theta_{k}}~=~\bigl(p^\theta_z(t, \bar z)\bfv\bigr)^{\dagger}\cdot \big[x^\theta_{zz}(t,\bar z)({\bf v}\otimes{\bf v})\big]_{\theta_{k}}\,,\\[4mm]
\big[\bigl(x^\theta_z(t, \bar z)\bfv\bigr)^{\dagger}\cdot p^\theta_{zz}(t,\bar z)({\bf v}\otimes{\bf v})\big]_{\theta_{k}}~=~\bigl(x^\theta_z(t, \bar z)\bfv\bigr)^{\dagger}\cdot \big[p^\theta_{zz}(t,\bar z)({\bf v}\otimes{\bf v})\big]_{\theta_{k}}\,,
\end{cases}
\]
for all $ t\in [0,T]$. Using \eqref{Ham2} and \eqref{Ham7} we thus obtain
\[\begin{split}
\dfrac {d}{dt}\big[\bigl(p^\theta_z&(t, \bar z)\bfv\bigr)^{\dagger}\cdot x^\theta_{zz}(t,\bar z)({\bf v}\otimes{\bf v})\big]_{\theta_{k}}~=~{d\over dt}\left[\bigl(p^\theta_z(t, \bar z)\bfv\bigr)^{\dagger}\cdot \big[x^\theta_{zz}(t,\bar z)({\bf v}\otimes{\bf v})\big]_{\theta_{k}}\right]\\ 
&=~-\bigl(x^\theta_z(t, \bar z)\bfv\bigr)^{\dagger}\cdot A(t)\cdot \big[x^\theta_{zz}(t,\bar z)({\bf v}\otimes{\bf v})\big]_{\theta_{k}}+\bigl(p^\theta_z(t, \bar z)\bfv\bigr)^{\dagger}\cdot C(t)\cdot\big[p^\theta_{zz}(t,\bar z)({\bf v}\otimes{\bf v})\big]_{\theta_{k}}\\
&={d\over dt}\left[\bigl(x^\theta_z(t, \bar z)\bfv\bigr)^{\dagger}\cdot \big[p^\theta_{zz}(t,\bar z)({\bf v}\otimes{\bf v})\big]_{\theta_{k}}\right].
\end{split}\]
Moreover,  using  \eqref{T-2} we obtain
\[\begin{split}&\big[\bigl(p^\theta_z(0, \bar z)\bfv\bigr)^{\dagger}\cdot x^\theta_{zz}(0,\bar z)({\bf v}\otimes{\bf v})\big]_{\theta_{k}}~=~\big[\bigl(p^\theta_z(T, \bar z)\bfv\bigr)^{\dagger}\cdot x^\theta_{zz}(T,\bar z)({\bf v}\otimes{\bf v})\big]_{\theta_{k}}\\
&\qquad -\big[\bigl(x^\theta_z(T, \bar z)\bfv\bigr)^{\dagger}\cdot p^\theta_{zz}(T,\bar z)({\bf v}\otimes{\bf v})\big]_{\theta_{k}}+\big[\bigl(x^\theta_z(0, \bar z)\bfv\bigr)^{\dagger}\cdot p^\theta_{zz}(0,\bar z)({\bf v}\otimes{\bf v})\big]_{\theta_{k}}=-6\bfv_k^3\,.
\end{split}\]
By   (\ref{r-S}) and (\ref{S-k}), we conclude that  the Jacobian matrix $D_\theta \Phi^{\psi^\theta}(\bar z, \bfv, 0)$ of partial derivatives w.r.t.~$\theta_{ij}, \theta_k$ contains $(n + 1)$ 
columns which form the $(n+1)\times (n+1)$ submatrix  
\[
\Lambda_{k}~=~\begin{bmatrix}\Tilde{S}_{i_0,k}\\[3mm] b_k\end{bmatrix},\qquad\qquad \qquad b_k~=~ \big(*,\cdots, *,~-6\bfv_k^3\big),
\]
such that 
\[
\det(\Lambda_k)~=~-6\bfv_k^3\cdot \det(S_{i_0}(0))\qquad\forall k\in \{1,\cdots, n\}.
\]
Since ${\bf v}\in S^{n-1}$, one has $\bfv_k\not= 0$ for some $k\in\{1,\cdots,n\}$.  Therefore $\mathrm{rank}(\Lambda_k)=n+1$.
In turn, this yields 
\bel{rank}
\mathrm{rank} \bigl[D_\theta \Phi^{\psi^\theta}(\bar z, \bfv, 0)\bigr]~=~n+1,
\eeq
as claimed.
\v
{\bf 6.} By continuity, there exists a  neighborhood ${\cal N}_{\bar z, \bfv}$ of $(\bar z, \bfv)$ such that
\bel{fullrank}
 \mathrm{rank} \bigl[ D_\theta  \Phi^{\psi^\theta}(z',\bfv',0)\bigr]\,=\,n+1\qquad\forall (z', \bfv')\in {\cal N}_{\bar z, \bfv}.
\eeq
Covering the compact set $\overline{B}_k\times S^{n-1}$ 
with finitely many open neighborhoods ${\cal N}_\ell=
{\cal N}_{\bar z^\ell, \bfv^\ell}$, $\ell=1,\ldots,N$, we consider the family of combined perturbations
\bel{psiall}\psi^\theta(z)~=~\psi(z) +\sum_{\ell=1}^N\eta\big(z-\bar{z}^{\ell}\big)\cdot\vp\big(\theta^{\ell},z-\bar{z}^{\ell}\big).\eeq
Here $\theta^{\ell}\doteq \bigl(\theta^\ell_{11}\cdots,\theta^\ell_{nn} , \theta^\ell_1,\cdots,  \theta^\ell_n\bigr)$  for every $\ell\in  \{1,\cdots, N\}$.  By the previous analysis,
(\ref{fullrank}) yields 
\[
\mathrm{rank}\bigl[D_{\theta^{\ell}}\Phi^{\psi^{\theta}}(z,{\bf v},0)\bigr]\,=\,n+1,\qquad (z,{\bf v})\in \mathcal{N}_{\ell}. 
\]
Hence,  the Jacobian matrix $D_{\theta}\Phi^{\psi^{\theta}}(z,{\bf v},0)$ of partial derivatives
w.r.t.~all combined variables $\theta = (\theta^\ell_{ij}, \theta^\ell_k)$ satisfies 
\[
\mathrm{rank}\bigl[ D_{\theta}\Phi^{\psi^{\theta}}(z,{\bf v},0)\bigr]\,=\,n+1,\qquad \forall~(z,{\bf v})\in B_k\times S^{n-1}.
\]
%
%
Again by continuity, we still have
\bel{fullrank2}
\mathrm{rank} \bigl[D_\theta \Phi^{\psi^\theta}(\bar z,\bfv,\theta)\bigr]~=~n+1
\eeq
for all $(\bar z,\bfv,\theta)\in \overline{B}_k\times S^{n-1}\times\R^{N(n^2+n)}$ with $|\theta|<\ve$
sufficiently small.   By the transversality theorem \cite{Bloom, GG},
this implies that, for a dense set of values $\theta$, the smooth map $\Phi^{\psi^\theta}$
at (\ref{Phi-psi})
is transversal to $\{\bf 0\}$, restricted to the domain $(\bar z, \bfv)\in \overline{B}_k\times S^{n-1}$. We conclude that the $\mathcal{M}_k\subset\C^4(\R^n)$ is dense.

\v

{\bf 8.} Finally, we remove the assumption of the global existence of solutions to  the Hamiltonian system (\ref{PMP1})-(\ref{PMP3}). For every $\nu\geq 1$,  consider  the system
\bel{PMP1v}
\dot x\,=\,\vp_{\nu}(x,p)\,H_p(x,p),\qquad\qquad 
\dot p\,=\,-\vp_{\nu}(x,p)\,H_x(x,p),
\eeq
where $\vp_{\nu}:\R^n\times\R^n\mapsto [0,1]$ is a smooth cutoff function satisfying
$$\left\{\bega{rl} \vp_{\nu}(x,p)=1\quad &\hbox{if} ~~|x|+|p|\leq \nu,\\[1mm]
\vp_{\nu}(x,p)=0\quad &\hbox{if} ~~|x|+|p|\geq \nu+1. \enda\right.$$
For every $z\in \R^n$, (\ref{PMP1v}),(\ref{PMP3}) admits a global solution $t\mapsto (x^{\nu}(t,z),p^{\nu}(t,z))$.  Define the $\C^1$ map  $\Phi^{\psi}:\R^d\times S^{n-1}\mapsto\R^n\times\R$ 
by 
\bel{Phi-psi-nu}
\Phi^{\psi}_{\nu}(z,{\bf v})~\doteq~\left(x^{\nu}_z(0,z){\bf v} \,, ~ \bigl(p^{\nu}_z(0, z)\bfv\bigr)^{\dagger}\cdot x^{\nu}_{zz}(0,z)({\bf v}\otimes{\bf v})\right).
\eeq
From the previous step, for every $k\geq 1$ the set 
\bel{M-k-nu}
\mathcal{M}^{\nu}_{k}~\doteq~\left\{\psi\in \mathcal{C}^4(\R^n) \,;~~ \Phi^{\psi}_{\nu}\Big|_{\overline{B}_k\times S^{n-1}}\text{ is transversal to } \{\bf0\}\right\}
\eeq
 is an open and dense subset of $\mathcal{C}^4(\R^n)$. Hence,  the set 
$$\M~\doteq~\bigcap_{\nu,k\geq 1} \M^{\nu}_k$$
is a countable intersection of open dense subsets of  $\C^4(\R^n)$. If $\psi\in \M$, then  for every   $k\geq 1$ the restriction of $\Phi^\psi_{\nu}$ to the domain $\overline{B}_k\times S^{n-1}$
is transversal to $\{{\bf 0}\}$. Consequently, the restriction of the set
\[
\Omega^{\nu}_{\psi}~\doteq~\Big\{(z,\bfv)\in\R^n\times S^{n-1}\,;~x^{\nu}_z(0, z)\bfv =0, \quad \bigl(p^{\nu}_z(0, z)\bfv\bigr)^{\dagger}\cdot x^{\nu}_{zz}(0,z)({\bf v}\otimes{\bf v})=0\Big\}
\] 
to  a neighborhood of $ \overline{B}_k\times S^{n-1}$ is an embedded manifold of dimension $n-2$. Since  $k\geq 1$ is arbitrary, $\Omega^{\nu}_{\psi}$   is    an embedded manifold  of dimension $n-2$. Next, let $(z,\bfv)\in \Omega_{\psi}$.   Then, by Remark~\ref{open}, there exists
an open neighborhood ${\cal N}$ of $(z,\bfv)$ such that 
\[
\Omega_\psi\cap {\cal N} ~=~\Omega^{\nu}_{\psi}\cap {\cal N}
\]
for all $\nu$ sufficiently large. This completes the proof.
\endproof

\section{Generic structure  of the optimal feedback control}
\label{sec:4}
\setcounter{equation}{0}

Consider again the optimal control problem (\ref{oc1})-(\ref{oc2}).
For every terminal cost $\psi\in \mathcal{C}^4(\R)$, we shall denote by  
\begin{itemize}
\item $\Gamma_\psi\subset\R^n$ the set of all conjugate points
$y = x(0, \bar z)$,  such that $t\mapsto x(t,\bar z)$  is a globally optimal trajectory.
\item $\mathcal{V}_{\psi}\subset\R^n$ the set of  all 
$y\in\R^n$ such that the optimization problem (\ref{oc1})-(\ref{oc2}) 
with initial data
$x(0)=y$
admits two or more globally  optimal solutions.
\end{itemize}

\begin{theorem}\label{t:41}  Let $f_0,\ldots, f_m$ and $L$ be given, satisfying 
{\bf (A1)}-{\bf (A2)}. 
 Then there exists a $\mathcal{G}_{\delta}$ subset $\mathcal{M}\subseteq\mathcal{C}^4(\R^n)$ with the following property:
\begin{itemize}
\item [(i)] The set $\Gamma_\psi$ is closed, with  locally bounded $(n-2)$-dimensional Hausdorff measure.
\item[(ii)]
Restricted to the complement $\R^n\setminus  \Gamma_\psi$, the set $\mathcal{V}_{\psi}$ of initial points  where
the optimization problem has multiple global solutions is locally contained in the union of finitely many 
embedded manifolds of dimension $n-1$.   
\end{itemize}

\end{theorem}

{\bf Proof.} {\bf 1.} To prove that $\Gamma_\psi$ is closed, consider a sequence of points $ y_k \in \Gamma_\psi$ with $y_k\to y$ as $k\to\infty$.    Assume that  $y_k=x(0,z_k)$, where $t\mapsto x(t, z_k)$ 
 is a globally optimal trajectory, and 
\[
\det \bigl[x_z(0,z_k)\bigr]~=~0,\qquad\quad\forall k\geq 1.\]
Thanks  to the uniform
bounds on optimal controls and optimal trajectories proved in Lemma~\ref{l:21}, by possibly extracting a subsequence 
we can assume the convergence 
$z_k\to\bar{z}$. By continuity, $t\mapsto x(t,\bar{z})$  is a globally optimal trajectory and 
\[
x(0, \bar z)\,=\,\lim_{k\to\infty}x(0, z_k)\,=\,\lim_{k\to\infty}y_k\,=\,y.
\]
Moreover, by Remark~\ref{open}, the map $t\mapsto x(t,z)$ is well-defined in $[0,T]$ for all $z$ in a neighborhood of $\bar{z}$. Hence, 
\[
\det\bigl[x_z(0,\bar z)\bigr]\,=\,\lim_{k\to\infty} \det \bigl[x_z(0,z_k)\bigr]\,=\,0,
\]
which implies $y\in \Gamma_\psi $ and the set $\Gamma_\psi$ is closed.
\v
{\bf 2.}
Next, by Theorem \ref{t:31} there exists a $\G_\delta$ set $\M\subset\C^4(\R^n)$ such that, for $\psi\in \M$, the set 
 $\Omega_\psi$ in (\ref{J}) is an  embedded manifold  of dimension $n-2$. For a given $\psi\in \M$,  
 we observe by \cite[Theorem 2.2]{BMN} that the set of all conjugate points satisfies
the inclusion
\[
\Gamma_{\psi}~\subseteq~ \Big\{x\bigl(0,z)\,;~~(z,{\bf v})\in\Omega_\psi~~\hbox{for some}~z\in \R^n,~\bfv\in S^{n-1}\Big\}.
\]
By Lemma~\ref{l:21} 
there exists a continuous function $\beta:\R_+\mapsto \R_+$ such that 
\[
|z|~\leq~\beta \bigl(|x(0,z)|\bigr),
\]
whenever the trajectory $t\mapsto x(t,z)$ is globally optimal for the problem (\ref{oc1})--(\ref{oc3}) with initial data
$y= x(0,z)$. 
Therefore, for every $r>0$ one has
\[
\Gamma_{\psi}\cap \overline{B}_r~\subseteq~\Big\{x(0,z) \,;~~(z,{\bf v})\in\Omega_\psi \quad 
\hbox{for some}~z\in \overline{B}_{\beta(r)}\,,~\bfv\in S^{n-1} \Big\}.
\]
Since the map $(z,\bfv)\mapsto x(0,z)$ is locally Lipschitz and $\Omega_\psi$ is an embedded manifold of dimension $n-2$, by the properties of Hausdorff measures \cite{EG} we conclude  that   the set  $\Gamma_\psi$  has locally bounded $(n-2)$-dimensional Hausdorff measure.  In particular, it is nowhere dense.
\v
{\bf 3.} For every $\delta>0$, we denote by $\mathcal{V}_{\psi}^{\delta}$ the set of points $y\in\R^n\backslash  \Gamma_\psi$ such that the optimization problem
has two distinct global solutions $t\mapsto x(t, z_i)$, $i=1,2$, with  
$$y\,=\,x(0, z_1)\,=\,x(0, z_2),\qquad\quad |z_1-z_2|~\geq~\delta.$$
Note that, by the previous steps, $\R^n\backslash  \Gamma_\psi$ is an open, everywhere dense set.

We claim that, for any compact set $K\subseteq \R^n\setminus  \Gamma_\psi$, there exists a constant $\delta_K>0$
such that
\[
 K\cap \mathcal{V}_{\psi}~\subseteq~ \mathcal{V}_{\psi}^{\delta_K}.
\]
Indeed, assume on the contrary that there exists a sequence of initial points  $y_k\in  K\cap \V_{\psi}$    from where two globally optimal trajectories start: 
 $t\mapsto x(t, z_{k,i})$, $i=1,2$, with  
$$y_k\,=\,x(0, z_{k,1})\,=\,x(0, z_{k,2}),\qquad |z_{k,1}-z_{k,2}|\to 0.$$
Setting $v_k=\ds {z_{k,2}-z_{k,1}\over |z_{k,2}-z_{k,1}|}$ and by possibly extracting a subsequence, thanks to Lemma \ref{l:21}, we can assume the convergence 
$$v_k\to \bar v,\qquad  y_k\to \bar y,\qquad z_{k,1}\,,z_{k,2}\to \bar z,$$
with  $|\bar{v}|=1$, and  $\bar y\in K$.  By continuity, $t\mapsto x(t,\bar{z})$  is a globally optimal trajectory.
Again, by Remark~\ref{open}, the map $t\mapsto x(t,z)$ is well-defined in $[0,T]$ for all $z$ in a neighborhood of $\bar{z}$. Thus, for $k\geq1$ large, we have  
\bel{v-D}\int_0^1 x_z\bigl( 0, \, \theta z_{k,2}+(1-\theta)z_{k,1}\bigr)(v_k) \, d\theta~=~\dfrac{x(0, z_{k,2})-x(0, z_{k,1})}{|z_{k,2}-z_{k,1}|}~=~0,
\eeq
and it follows
\[
x_z(0,\bar{z})(\bar v) ~=~\lim_{k\to\infty}x_z(0,z_{k,1}) (v_{k})~=~0.
\]
This implies $\bar{y}\in \Gamma_{\psi}$, yielding a contradiction. 
\v
{\bf 4.} To prove (ii), it remains to show that $K\cap \mathcal{V}_{\psi}^{\delta_K}$  is  contained in the union of finitely many embedded manifolds of dimension $n-1$. For every point $z\in\R^n$ which is the terminal point of a globally optimal trajectory $t\mapsto x(t,z)$ of (\ref{PMP1})--(\ref{PMP3}), call  $W(z)$ the cost of this trajectory. Namely,
\bel{W}
W(z)~=~\int_0^{T}L\Big(x(t,z),\,u\bigl(x(t,z),p(t,z)\bigr)\Big)dt+\psi(z).
\eeq
By Lemma \ref{l:21} and the compactness of $K$, there exists a constant $r_K>0$ sufficiently large such that 
\[\bega{l} 
K\cap \mathcal{V}_{\psi}^{\delta_K}~=~
\Big\{y\in K\,;~~\hbox{there exist $z_1, z_2$ such that}~y=x(0,z_1)=x(0, z_2), ~~ |z_1-z_2|\geq\delta_K,\\[3mm]
\qquad\quad \hbox{moreover, for $i=1,2$,} ~|z_i|\leq r_K ~ \hbox{and the trajectory $t\mapsto  x(t, z_i)$ is globally optimal}\Big\}.
\enda
\]
\v
{\bf 5.} We now observe that compact set of couples
$$\Big\{ (z_1, z_2)\in \R^n\times\R^n\,;~|z_1|\leq r_K\,,~|z_2|\leq r_K \,,~~|z_1-z_2|\geq\delta_K\Big\}$$
can be covered by finitely many open sets of the form $B(\bar z_1, \delta_K/2)\times B(\bar z_2, \delta_K/2)$,
where $B(z,\delta)$ denotes the open ball centered at $z$ with radius $\delta$.

To complete the proof of (ii), it thus suffices to prove the following statement.
\begi
\item[{\bf (P)}] {\it For every $(\bar{z}_1,\bar{z}_2)\in \R^n\times\R^n$ with $|\bar{z}_1-\bar{z}_2|\geq \delta_K>0$, the set of initial points $y\in K$
such that 
$$y\,=\,x(0, z_1)\, =\, x(0, z_2), \qquad\qquad W(z_1)~=~W(z_2),$$
for some $z_i\in B(\bar z_i, \delta_K/2)$, $i=1,2$, is contained in an embedded manifold of 
dimension $n-1$.}
\end{itemize}
Indeed, $(z_1,z_2)\in B(\bar{z}_1,\delta_K/2)\times B(\bar{z}_2,\delta_K/2)$ implies  $z_1\neq z_2$. Hence, by the uniqueness of solutions to (\ref{PMP1}), if $x(0, z_1)=x(0, z_2)$,  we must have $p(0, z_1)\not= p(0, z_2)$. 
This implies
\bel{rk} \nabla W(z_1) 
\bigl[x_z(0, z_1)\bigr]^{-1} ~= ~ p(0,z_1)~\not= ~p(0,z_2) ~= ~ \nabla W(z_2) \bigl[x_z (0, z_2)\bigr]^{-1}.
\eeq
Consider  the $\C^1$ map $\Phi:\R^{n}\times\R^n\mapsto\R^n\times\R$ 
defined by
\bel{Phi}
\Phi(z_1,z_2)~\doteq~\Big(x(0,z_1)-x(0,z_2) \,, ~ W(z_1)-W(z_2)\Big).
\eeq
We  claim that the $(n+1)\times 2n$ Jacobian matrix
\[D\Phi(z_1,z_2)~=~\begin{bmatrix} x_z(0,z_1)& -x_z(0,z_2)\\[4mm] \nabla W(z_1)&-\nabla W(z_2)
\end{bmatrix}
\]
has maximum rank. Otherwise,
 its last row would be a linear combination of the other rows, i.e. 
 there would exist ${\bf u}\in\R^n$ such that
$$\left\{\begin{array}{ll}\nabla W(z_1)~=~{\bf u}\cdot x_z(0,z_1),\\[2mm]
 \nabla W(z_2)~=~{\bf u}\cdot x_z(0,z_2),\end{array}\right.$$
in contradiction with (\ref{rk}). Therefore, $\rank\bigl[D\Phi(z_1,z_2)\bigr]=n+1$, and {\bf (P)} follows by the 
implicit function theorem \cite{Bloom, GG}.
\endproof
\medskip

{\small
{\bf Acknowledgment.}  The research of K.\,T.\,Nguyen was partially supported by NSF with grant DMS-2154201, ``Generic singularities and fine regularity structure for nonlinear partial differential equations".}

\end{document}